\documentclass[a4paper]{amsart}
\usepackage{amsmath,amsfonts,amssymb,amsthm}
\usepackage{mathtools}
\usepackage{mathrsfs}
\usepackage{pdfpages}
\usepackage{url}
\usepackage{color,xcolor}
\usepackage{caption}
\usepackage{enumerate}
\usepackage{afterpage}
\usepackage{xargs,twoopt}
\usepackage{hyperref}
\usepackage[toc,page]{appendix}
\usepackage{spverbatim}
\hypersetup{colorlinks = true, linkbordercolor = {white}}
\usepackage{aliascnt}
\usepackage{bbold}
\usepackage{dsfont}
\usepackage[cal=boondoxo]{mathalfa}
\usepackage{bussproofs}
\usepackage{graphicx}
\usepackage[retainorgcmds]{IEEEtrantools}
\usepackage{stmaryrd}
\usepackage{tensor}
\usepackage{textcomp}
\usepackage{tikz}
\usepackage{units}
\usepackage[all]{xy}
\usepackage{comment}

\DeclareSymbolFont{extraitalic}      {U}{zavm}{m}{it}
\DeclareMathSymbol{\Qoppa}{\mathord}{extraitalic}{161}
\DeclareMathSymbol{\qoppa}{\mathord}{extraitalic}{162}
\DeclareMathSymbol{\Stigma}{\mathord}{extraitalic}{167}
\DeclareMathSymbol{\Sampi}{\mathord}{extraitalic}{165}
\DeclareMathSymbol{\sampi}{\mathord}{extraitalic}{166}
\DeclareMathSymbol{\stigma}{\mathord}{extraitalic}{168}

\usepackage{tikz,tikz-cd}
\usetikzlibrary{matrix}
\usetikzlibrary{fit}
\usetikzlibrary{calc}
\usetikzlibrary{decorations.pathmorphing}
\tikzset{curve/.style={settings={#1},to path={(\tikztostart)
    .. controls ($(\tikztostart)!\pv{pos}!(\tikztotarget)!\pv{height}!270:(\tikztotarget)$)
    and ($(\tikztostart)!1-\pv{pos}!(\tikztotarget)!\pv{height}!270:(\tikztotarget)$)
    .. (\tikztotarget)\tikztonodes}},
    settings/.code={\tikzset{quiver/.cd,#1}
        \def\pv##1{\pgfkeysvalueof{/tikz/quiver/##1}}},
    quiver/.cd,pos/.initial=0.35,height/.initial=0}

\tikzset{tail reversed/.code={\pgfsetarrowsstart{tikzcd to}}}
\tikzset{2tail/.code={\pgfsetarrowsstart{Implies[reversed]}}}
\tikzset{2tail reversed/.code={\pgfsetarrowsstart{Implies}}}
\tikzset{no body/.style={/tikz/dash pattern=on 0 off 1mm}}

\newcommand{\circv}{\circ_{\mathsf{v}}}
\newcommand{\circh}{\circ_{\mathsf{h}}}

\newcommand{\ceJ}{\mathcal{J}}

\newcommand{\ceL}{\mathcal{L}}

\newcommand{\cB}{\mathscr{B}}
\newcommand{\cC}{\mathscr{C}}

\newcommand{\cI}{\mathscr{I}}

\newcommand{\cV}{\mathscr{V}}

\newcommand{\mbC}{\mathbf{C}}

\newcommand{\mbM}{\mathbf{M}}
\newcommand{\mbN}{\mathbf{N}}

\newcommand{\bbon}{\mathbb{1}}
\newcommand{\bbk}{\mathbb{k}}

\newcommand{\tti}{\mathtt{i}}
\newcommand{\ttj}{\mathtt{j}}

\newcommand{\ol}[1]{\overline{#1}}
\newcommand{\on}[1]{\operatorname{#1}}

\DeclareMathOperator{\End}{End}

\DeclareMathOperator{\Hom}{Hom}
\DeclareMathOperator{\HH}{HH}
\DeclareMathOperator{\Id}{Id}

\newcommand{\Vect}{{\cV}\hspace{-2.0pt}ec}

\DeclareMathOperator{\add}{add}
\DeclareMathOperator{\afmod}{\text{-}afmod}

\DeclareMathOperator{\ev}{ev}

\DeclareMathOperator{\im}{im}

\DeclareMathOperator{\id}{\mathsf{id}}

\DeclareMathOperator{\modd}{mod}
\DeclareMathOperator{\bmodd}{\textbf{mod}}

\DeclareMathOperator{\proj}{proj}
\DeclareMathOperator{\bproj}{\textbf{proj}}

\EnableBpAbbreviations

\newtheorem{thm}{Theorem}[section]
\newaliascnt{lem}{thm}
\newtheorem{lem}[lem]{Lemma}
\aliascntresetthe{lem}
\newaliascnt{psn}{thm}
\newtheorem{psn}[psn]{Proposition}
\aliascntresetthe{psn}
\newaliascnt{cor}{thm}
\newtheorem{cor}[cor]{Corollary}
\aliascntresetthe{cor}
\newaliascnt{que}{thm}

\aliascntresetthe{que}
\newaliascnt{conj}{thm}

\aliascntresetthe{conj}
\theoremstyle{definition}
\newaliascnt{ex}{thm}

\aliascntresetthe{ex}
\newaliascnt{exs}{thm}

\aliascntresetthe{ex}
\newaliascnt{def}{thm}
\newtheorem{defn}[def]{Definition}
\aliascntresetthe{def}
\newaliascnt{nota}{thm}

\aliascntresetthe{nota}
\newaliascnt{rmk}{thm}
\newtheorem{rmk}[rmk]{Remark}
\aliascntresetthe{rmk}
\newaliascnt{rmks}{thm}

\aliascntresetthe{rmks}

\DeclareMathOperator{\Ext}{Ext}

\usepackage{stmaryrd}
\title{Hochschild cohomology for finitary $2$-representations}
\author{James Macpherson, Vanessa Miemietz, Mateusz Stroi\'nski }
\begin{document}

\maketitle

\begin{abstract}
In this article, we define and investigate Hochschild cohomology for finitary $2$-representations of quasi-fiat $2$-categories.
\end{abstract}

\section*{Introduction}

Since the beginning of this century, major progress in representation theory has been obtained by means of categorification. Loosely speaking, one replaces actions of algebraic structures on vector spaces via linear transformations by actions on categories via functors. These functors have natural transformations between them, providing an extra level of structure. Categorifications of quantum groups \cite{CR, KL, Ro} have led to the proof of Brou\'e's abelian defect group conjecture for symmetric groups, as well as major progress in the representation theory of Hecke algebras. Meanwhile, categorifications of Hecke algebras \cite{Soe} have led among other progress to the proof of Kazhdan--Lusztig conjectures for arbitrary Coxeter groups \cite{EW} and counter-examples to Lusztig's and James' conjectures.

The success of categorification has led to the development of approaches to the representation theory of suitable kinds of $2$-categories or bicategories, e.g. those of finitary or fiat $2$-categories initiated in \cite{MM1} or those of tensor categories in \cite{EGNO}. In both cases, one of the stepping stones is the observation that $2$-representations can (under mild assumptions) be internalised and realised as certain categories of (co)modules over (co)algebras \cite{MMMT}.

An important homological invariant for algebras is given by Hochschild cohomology \cite{Ho}. The Hochschild cohomology groups form a graded ring under cup product, as well as a Lie algebra under the Gerstanhaber bracket. In degree zero, it describes the centre of the algebra; in degree one, it is given by the quotient of the space of derivations by inner derivations; and provided the algebra is defined over a field, in higher degrees it determines the deformation theory of the algebra. This comes from the fact that, for an algebra $A$ over a field, the classical definition via the bar resolution can be replaced by the description of Hochschild cohomology as the Yoneda extension algebra of $A$ in the category of $A$-$A$-bimodules. 

General Hochschild cohomology has been studied for algebra or ring objects in monoidal categories in \cite{HF1,HF2}, by generalising the bar resolution. In this article, we take a slightly different approach, focusing on algebra $1$-morphisms in (quasi-)fiat $2$-categories and using the definition of Hochschild cohomology as the Yoneda extension algebra of $A$ in the category of $A$-$A$-bimodules.  In this way, Hochschild cohomology becomes an invariant of the $2$-representation determined by a given algebra $1$-morphism (see Proposition \ref{moritainvariant}).
Both approaches coincide if the corresponding $2$-representation is simple transitive in the sense of \cite{MM5} (hereafter, we call this \emph{simple}). In general, an algebra $1$-morphism associated to a finitary $2$-representation of a quasi-fiat $2$-category $\cC$ does not live in $\cC$ itself, but rather in an abelianisation $\ol{\cC}$. In this abeliansation, the projective $1$-morphisms are precisely those coming from $\cC$. In particular, if an algebra $1$-morphism does not belong to the non-abelanised $2$-category $\cC$, the bar resolution does not form a projective resolution of $A$, leading us to investigate a replacement of the bar resolution.

 For any finitary $2$-representation, zeroth and first Hochschild cohomology provide the natural analogues of centres and derivations of the corresponding algebras (see Propositions \ref{centre} and \ref{deriv}). As for algebras over a field, the centre of a simple algebra $1$-morphism is trivial. On the other hand, we provide an example of a simple algebra $1$-morphism with nontrivial first Hochschild cohomology in Section \ref{VecGHH}. In a departure from the classical picture, second Hochschild cohomology not only takes into account a generalisation of the usual cocycle condition for algebras over a field, but also has an ingredient coming from extensions of Yoneda degree $1$ of the algebra with itself in the ambient (abelianised) $2$-category, see Proposition \ref{HH2}. Moreover, second Hochschild cohomology only gives rise to deformations under additional assumptions (Proposition \ref{HH2def}). 

The article is structured as follows: In Section \ref{bckg}, we recall the necessary background from finitary $2$-representation theory. In Section \ref{mainsec}, we define Hochschild cohomology, prove Morita invariance and analyse Hochschild cohomology groups in small degrees. In Section \ref{examplesec}, we provide examples.

\section{Background}\label{bckg}

Throughout this article, let $\Bbbk$ be an algebraically closed field.

\subsection{Multifinitary $2$-categories and finitary $2$-representations.}

A \emph{finitary category} is 
an additive $\Bbbk$-linear category, which is idempotent complete, has only finitely many indecomposable objects up to isomorphism, and whose morphism spaces are finite-dimensional. Equivalently, it is the Cauchy completion of a one-object category given by a finite-dimensional algebra.

We denote by $\mathfrak{A}_{\Bbbk}^{f}$ the $2$-category whose objects are finitary categories, whose $1$-morphisms are $\Bbbk$-linear functors and whose $2$-morphisms are natural transformations of such functors. 

A \emph{multifinitary $2$-category} is a $2$-category with finitely many objects, in which all morphism categories are finitary and horizontal composition of $2$-morphisms is $\Bbbk$-bilinear.

We will usually denote objects in multifinitary $2$-categories by $\bullet$ (if there is only one) or $\tti, \ttj,$ etc., $1$-morphisms by $F,G,$ etc., and $2$-morphisms by $\alpha, \beta,$ etc..

A multifinitary $2$-category $\cC$ is called \emph{quasimultifiat} if there exists a biequivalence ${}^*\colon\cC\to \cC^{co,op}$, such that for any $1$-morphism $F\in \cC(\tti,\ttj)$, there is a  unit $\bbon_{\tti} \to F^*F$ and a counit $FF^*\to\bbon_{\ttj}$ satisfying the usual adjunction identity. A quasimultifiat $2$-category is \emph{multifiat} if ${}^*$ is weakly involutive.

A \emph{finitary $2$-representation} $\mbM$ of a multifinitary $2$-category $\cC$ is a $2$-functor $\mbM\colon \cC\to \mathfrak{A}_\Bbbk^f$. Finitary $2$-representations of $\cC$ together with morphisms of $2$-representations (strong natural transformations of $2$-functors) and modifications form a $2$-category, denoted by $\cC\afmod$.

A finitary $2$-representation $\mbM$ of a multifinitary $2$-category $\cC$ is called \emph{cyclic} if there exists an object $X\in \mbM(\tti)$ for some $\tti\in\cC$, such that $\add \{\mbM(F)\,X\,\vert\, F\in \cC(\tti,\ttj)\}\simeq \mbM(\ttj)$ for all $\ttj\in \cC$. Such an $X$ is called a \emph{generator} of $\mbM$.  The $2$-representation $\mbM$ is called \emph{transitive} if it is generated by any $X$ in any of the $\mbM(\tti)$.

A finitary $2$-representation $\mbM$ of a multifinitary $2$-category $\cC$ is called \emph{simple} if it has no nontrivial ideals, where an \emph{ideal} consists of a collection of categorical ideals in each $\mbM(\tti)$, which is stable under the action by $\cC$.

\subsection{Cells.} One of the most powerful tools in order to classify simple $2$-representations of a given multifinitary $2$-category $\cC$ is the use of cell structures. Given two indecomposable $1$-morphisms $F$ and $G$, we say $F\leq_LG$ if there exists a $1$-morphism $H$ such that $G$ is a direct summand of $HF$. We call the resulting partial preorder the \emph{left order} and the corresponding equvialence classes \emph{left cells}. Analogously, we define the \emph{right order} by saying $F\leq_RG$ if there exists an $H$ such that $G$ is a direct summand of $FH$, and the corresponding \emph{right cells}, as well as the \emph{two-sided order} where $F\leq_JG$ if there exist $H_1,H_2$ such that $G$ is a direct summand of $H_1FH_2$, and the corresponding \emph{two-sided cells}. A two-sided cell $\ceJ$ is called \emph{strongly regular} if the intersection of any left an any right cell contained in it only contains one isomorphism class of indecomposable $1$-morphism.

Any transitive (and hence any simple) $2$-representation $\mbM$ of a multifinitary $2$-category $\cC$ has an \emph{apex}, which is the unique two-sided cell $\ceJ$ that is maximal with respect to the condition that $\mbM(\ceJ)\neq 0$. A $2$-category $\cC$ is called \emph{$\ceJ$-simple} for a two-sided cell $\ceJ$, if any nonzero $2$-ideal in $\cC$ necessarily contains the identity $2$-morphisms on the $1$-morphisms in $\ceJ$. By factoring out the maximal $2$-ideal not containing the identity $2$-morphisms on the $1$-morphisms in $\ceJ$, we obtain the so-called \emph{$\ceJ$-simple quotient}.

\subsection{Abelianisation.} In order to internalise $2$-representations, we need the concept of abelianisation. We denote by $\ol{\mathcal{C}}$ the \emph{projective abelianisation}, as introduced by Freyd \cite{Fr}, of a finitary category $\mathcal{C}$, whose objects are morphisms $X_1\xrightarrow{x}X_0$ in $\cC$, and whose morphisms are pairs $(f_0,f_1)$ given by (solid) commutative diagrams of the form
$$\xymatrix{
X_1\ar^{x}[rr]\ar^{f_1}[d]&&X_0\ar^{f_0}[d]\ar@{-->}_{h}[lld]\\
Y_1\ar^{y}[rr]&&Y_0
}$$
modulo the homotopy relation that such a pair $(f_0,f_1)$ defines the zero morphism if there exists an $h$ such that $f_0=yh$. 

Similarly, we can define the \emph{projective abelianisation} $\ol{\cC}$ of a multifinitary $2$-category $\cC$, by setting $\ol{\cC}(\tti,\ttj) = \ol{\cC(\tti,\ttj)}$. Horizontal composition is given by 
$$(F_1\xrightarrow{\alpha}F_0)(G_1\xrightarrow{\beta}G_0) = (F_1G_0\oplus F_0G_1\xrightarrow{(\alpha\circh\id_{G_0}, \id_{F_0}\circh \beta)}F_0G_0).$$

This technically only defines a bicategory and there is a more technical version producing a $2$-category, see \cite[Section 3.2]{MMMT} for the dual version, but the above definition will suffice for the purpose of this article.

Similary, given a finitary $2$-representation $\mbM$ of a multifinitary $2$-category $\cC$, we can define its abelianisation $\ol{\mbM}$ by $\ol{\mbM}(\tti) = \ol{\mbM(\tti)}$ with the action of $\cC$, or $\ol{\cC}$, given component-wise. 

\subsection{$2$-categories of $2$-representations.}

Finitary $2$-representations of a fixed $2$-category $\cC$ again form a $2$-category, whose $1$-morphisms are morphisms of $2$-representations (strong natural transformations) and whose $2$-morphisms are modifications. This $2$-category is denoted by $\cC\afmod$. 

Given two finitary $2$-representations $\mbM,\mbN$, the morphism category $\Hom_{\cC}(\mbM, \mbN)$ has a full subcategory $\Hom_{\cC}^{ex}(\mbM, \mbN)$ given by \emph{exact morphisms}, where a morphism is called exact if each component functor $\mbM(\tti)\to \mbN(\tti)$ induces an exact functor on the abelianisation.
The $2$-full sub-$2$-category on the same objects, but whose $1$-morphisms are only the exact morphisms is denoted by $\cC\afmod^{ex}$.

\subsection{Algebra $1$-morphisms and modules}

Throughout this section, let $\cC$ be a finitary $2$-category.

An \emph{algebra $1$-morphism} $(A,\mu,\iota)$ in $\ol{\cC}$ is a monoid object in the category $\ol{\cC}(\tti,\tti)$, for some $\tti \in \ol{\cC}$. In category-theoretic settings, one often refers to algebra $1$-morphisms in a $2$-category as {\it monads} therein.
 
Explicitly, it is a $1$-morphism $A\colon \tti\to\tti$ in $\ol{\cC}$ along with $2$-morphisms $\mu\colon A\circ A\to A$ and $\iota\colon \bbon_\tti\to A$ subject to the following standard algebra axioms:\begin{itemize}
		\item $\mu\circv(\id_A\circh\mu)=\mu\circv(\mu\circh\id_A)$;
		\item $\mu\circv(\id_A\circh\iota)=\mu\circv(\iota\circh\id_A)=\id_A$.
\end{itemize}

Let $\cC$ be a finitary $2$-category and let $A\colon \tti\to\tti$ be an algebra $1$-morphism in $\ol{\cC}$. A \emph{(right) $A$-module $1$-morphism} $(M,\rho_M)$ is a $1$-morphism $M\colon \tti\to\ttj$ and a $2$-morphism $\rho:=\rho_M\colon M\circ A\to M$ in $\ol{\cC}$ such that the following standard module axioms hold:\begin{itemize}
		\item $\rho\circv(\rho\circh\id_A)=\rho\circv(\id_M\circh\mu)$;
		\item $\rho\circv(\id_M\circh\iota)=\id_M$.
	\end{itemize}

A \emph{morphism of (right) $A$-module $1$-morphisms} $\alpha\colon (M,\rho_M)\to (N,\rho_n)$ is a $2$-morphism $\alpha\colon M\to N$ of $\ol{\cC}$ such that $\rho_N\circv(\alpha\circh\id_A)=\alpha \circv \rho_M$. We can similarly define \emph{left $A$-module $1$-morphisms} $(M,\lambda_M)$ of $A$, with $\lambda_M\colon  A\circ M\to M$ with analogous axioms and morphisms. We denote the category of right $A$-module $1$-morphisms in $\ol{\cC}(\tti,\ttj)$ by $\modd_{\ol{\cC}(\tti,\ttj)}\text{-}A$, and the category of left $A$-module $1$-morphisms in $\ol{\cC}(\ttj,\tti)$ by $A\text{-}\modd_{\ol{\cC}(\ttj,\tti)}$.

A \emph{bimodule $1$-morphism} $(M,\rho_M,\lambda_M)$ of $A$ is a $1$-morphism $M$ in $\ol{\cC}$ such that $(M,\rho_M)$ is a right $A$-module, $(M,\lambda_M)$ is a left $A$-module, and $\rho_M\circv(\lambda_M\circh\id_A)=\lambda_M\circv(\id_A\circh\rho_M)$. A \emph{morphism} of $A$-$A$-bimodules $\alpha\colon (M,\rho_M,\lambda_M)\to (N,\rho_N,\lambda_N)$ is a $2$-morphism $\alpha\colon M\to N$ in $\ol{\cC}$ that is both a left and a right $A$-module morphism.
We denote the category of $A$-$A$-bimodule $1$-morphisms in $\ol{\cC}(\tti,\tti)$ by $A$-$\modd_{\ol{\cC}(\tti,\tti)}$-$A$.

For algebra $1$-morphisms $A,B$ and $C$ and bimodules $(M,\rho_M,\lambda_M)\in A\text{-}\modd_{\ol{\cC}(\ttj,\tti)}\text{-}B$ and $(N,\rho_N,\lambda_N)\in B\text{-}\modd_{\ol{\cC}(\tti,\ttj)}\text{-}C$, we define $M\circ_B N$ as the cokernel of $$MBN\xrightarrow{\rho_M\circh\id_N-\id_M\circh\lambda N}MN. $$

We will often just talk about algebras, modules and bimodules, without always carrying around the word $1$-morphism. We will also often omit the explicit mention of the action(s) in the name of the (bi)module.

\subsection{Free-forgetful adjunctions}\label{sec:freeforgadj}
Fix an algebra $1$-morphism $A\in \ol{\cC}(\tti,\tti)$ and let $(M, \rho_M, \lambda_M)\in A$-$\modd_{\ol{\cC}(\tti,\tti)}$-$A$. The standard free-forgetful adjunctions give rise to isomorphisms
$$\begin{array}{rcccl}
\Hom_{\ol{\cC}(\tti,\tti)}(F,M)& \cong& \Hom_{\modd_{\ol{\cC}(\tti,\tti)}\text{-}A} (FA,M)& \cong& \Hom_{A\text{-}\modd_{\ol{\cC}(\tti,\ttj)}\text{-}A} (AFA,M)\\
f &\mapsto &\rho_M \circv(f\circh \id_A)&&\\
g\circv(\id_F\circh \iota_A)&\mapsfrom &g&&\\
&&g &\mapsto & \lambda_M\circv(\id_A\circh g)\\
&&h\circv(\iota_A\circh \id_{FA})&\mapsfrom &h\\
\end{array}$$ 
for $F\in \ol{\cC}(\tti,\tti)$.

\subsection{Categories of modules as $2$-representations}

The categories $\modd_{\ol{\cC}(\tti,\ttj)}$-$A$, $A$-$\modd_{\ol{\cC}(\ttj,\tti)}$ and $A$-$\modd_{\ol{\cC}(\tti,\tti)}$-$A$ are abelian categories, and we can consider their subcategories of projective objects. We notate these subcategories as $\proj_{\ol{\cC}(\tti,\ttj)}$-$A$, $A$-$\proj_{\ol{\cC}(\ttj,\tti)}$ and $A$-$\proj_{\ol{\cC}(\tti,\tti)}$-$A$, respectively.

For a quasimultifiat $2$-category $\cC$ and an algebra $1$-morphism $A$, there is a finitary $2$-representation
$\bproj_{\ol{\cC}}$-$A$ given by $\bproj_{\ol{\cC}}$-$A(\ttj) = \proj_{\ol{\cC}(\tti,\ttj)}$-$A$, where the $\cC$-action is just the natural action on the left.
The abelianisation of $\bproj_{\ol{\cC}}$-$A$ is $\bmodd_{\ol{\cC}}$-$A$, where $\bmodd_{\ol{\cC}}$-$A(\ttj) = \modd_{\ol{\cC}(\tti,\ttj)}$-$A$.

 For an algebra $1$-morphism $A\in \cC(\tti,\tti)$, the category $\proj_{\ol{\cC}(\tti,\ttj)}$-$A$ is given by the additive closure of $\add\{FA \mid \, F\in \cC(\tti,\ttj)\}$ inside $\modd_{\ol{\cC}(\tti,\ttj)}$-$A$, see e.g.\ \cite[Lemma 4.1]{MMMTZbicat}. Likewise, $A$-$\proj_{\ol{\cC}(\ttj,\tti)}$  is given by the additive closure inside $A$-$\modd_{\ol{\cC}(\ttj,\tti)}$ of $\add\{AF \mid \, F\in \cC(\ttj,\tti)\}$.

 Similarly, one shows the analogous results for bimodules.

\begin{lem}\label{lem:ProjBimodDefn}
The category $A$-$\proj_{\ol{\cC}(\tti,\tti)}$-$A$ is given by
 $\add\{AFA\mid \, F\in \cC(\tti,\tti)\}$ inside $A$-$\modd_{\ol{\cC}(\tti,\tti)}$-$A$.
\end{lem}

\proof
Using the free-forgetful adjunction, $$ \Hom_{A\text{-}\modd_{\ol{\cC}(\tti,\ttj)}\text{-}A} (AFA,-)\cong \Hom_{\ol{\cC}(\tti,\tti)}(F,-),$$ which is exact if $F\in \cC(\tti,\tti)$. The proof that every bimodule is indeed a quotient of some bimodule of the form $AFA$ for  $F\in \cC(\tti,\tti)$ is analogous to the proof for right comodules in  \cite[Lemma 4.1]{MMMTZbicat}.
\endproof

\subsection{Internal hom}\label{inthomsec}

For a quasi-fiat $2$-category $\cC$, a finitary $2$-representation $\mbM$ of $\cC$, and an object $X$ in one of the $\mbM(\tti)$, 
the evaluation functors  $\ol{\mbM(-)}X \colon \ol{\cC}(\tti,\ttj) \to \ol{\mbM}(\tti), \; F\mapsto \ol{\mbM}(F)X$ are right exact,
and their right adjoint is called internal hom and denoted by $[X,-]$. In particular, for all $F\in \ol{\cC}(\tti,\ttj), Y\in \ol{\mbM}(\ttj)$, we have an isomorphism
$$\Hom_{\ol{\cC}(\tti,\ttj)} (F, [X,Y]) \cong \Hom_{\ol{\mbM}(\ttj)}(\mbM(F)X,Y).$$

By the dual arguments from  \cite[Section~4]{MMMT} (cf.\ also \cite[Section~7.10]{EGNO}, \cite[Section~4]{LM}), the internal hom $[X,X]\in \ol{\cC}(\tti,\tti)$ carries the structure of an algebra $1$-morphism and, provided $X$ generates $\mbM$, the (collection of) functor(s) $[X,-]$ defines an equivalence of $2$-representations between $\mbM$ and $\bproj_{\ol{\cC}}$-$[X,X]$.

Note that, for a $2$-morphism $\alpha \colon F\to G$ in $\ol{\cC}(\tti,\ttj)$ and $Y\in \ol{\mbM}(\ttj)$ functoriality of the internal hom implies commutativity of the diagram
$$ \xymatrix{
\Hom_{\ol{\cC}(\tti,\ttj)}(G, [X,Y]) \ar@{<->}^{\sim}[rr] \ar^{-\circv \alpha}[d]
&&\Hom_{\ol{\mbM}(\ttj)}(\mbM(G)X,Y)  \ar^{-\circv \mbM(\alpha)_X}[d]\\
\Hom_{\ol{\cC}(\tti,\ttj)}(F, [X,Y]) \ar@{<->}^{\sim}[rr]
&& \Hom_{\ol{\mbM}(\ttj)}(\mbM(F)X,Y).
}$$

\section{Hochschild cohomology of algebra $1$-morphisms}

\subsection{Hochschild cohomology}\label{mainsec}

Let $\cC$ be a quasi-fiat $2$-category, and let $\mbM$ be a cyclic finitary $2$-representation of $\cC$. By Section \ref{inthomsec} above, there exists an algebra $1$-morphism $A=A_{\mbM}$ such that $\mbM$ is equivalent to $\bproj_{\ol{\cC}}$-$A$ as $2$-representations of $\cC$. 
Consequently, for the sequel, we will be 
commonly referring to (abelianisations of) $2$-representations as their equivalent internal module categories.

For the rest of the section, fix $\mbM$ and $A=A_\mbM$.\par

In this case, $A$-$\modd_{\ol{\cC}(\tti,\tti)}$-$A$ is an abelian category, and we can construct projective resolutions of bimodules, and more generally the derived category $\mathcal{D}(A\text{-}\modd_{\ol{\cC}(\tti,\tti)}\text{-}A)$.

\begin{defn} For $k\geq 0$, the \emph{$k$th Hochschild cohomology} of $A$ is is defined to be $$\HH_{\cC}^k(A) = \Hom_{\mathcal{D}(A\text{-}\modd_{\ol{\cC}(\tti,\tti)}\text{-}A)}(A, A[k]),$$ or, equivalently,
$$\HH_{\cC}^*(A) = \Ext^*_{A\text{-}\modd_{\ol{\cC}(\tti,\tti)}\text{-}A}(A, A).$$\end{defn}

More explicitly, let $\dots \to P_2\to P_1\to P_0\to A$ be a projective resolution of $A$ in $A\text{-}\modd_{\ol{\cC}(\tti,\tti)}\text{-}A$. Then $\HH_{\cC}^k(A)$ can be calculated as the $k$th cohomology of the complex $\Hom_{A\text{-}\modd_{\ol{\cC}(\tti,\tti)}\text{-}A}(P_{\bullet},A)$.

Note that while we have defined Hochschild cohomology for algebra $1$-morphisms, the following easy lemma shows that Hochschild cohomology is a Morita invariant and it thus makes sense to talk about the Hochschild cohomology of a given $2$-representation.

\begin{psn}\label{moritainvariant}
Let $\mbM$ be a finitary $2$-representation of $\cC$ and $A\in \ol{\cC}(\tti,\tti)$ and $B\in \ol{\cC}(\ttj,\ttj)$ two algebra $1$-morphisms  such that $\mbM_A\cong \mbM \cong \mbM_B$. Then $\HH_\cC^*(A)\cong \HH_\cC^*(B)$.
\end{psn}

\proof
By \cite[Theorem 19]{MMMT}, we may assume that there exist biprojective bimodules $M\in A\text{-}\modd_{\ol{\cC}(\ttj,\tti)}\text{-}B, N\in B\text{-}\modd_{\ol{\cC}(\tti,\ttj)}\text{-}A$ such that $M\circ_BN\cong A$ and $N\circ_A M\cong B$. 
In particular $(M\circ_B - , N\circ_A - )$ and $(-\circ_A M, -\circ_B N)$ form biadjoint equivalences. Then $M \circ_{B} - \circ_{B} N$ becomes an equivalence, with a quasi-inverse given by $N \circ_{A} - \circ_{A} M$. We find 
\[
\begin{aligned}
&\HH_{\cC}^k(A) = \Hom_{\mathcal{D}(A\text{-}\modd_{\ol{\cC}(\tti,\tti)}\text{-}A)}(A, A[k]) \\
&\simeq \Hom_{\mathcal{D}(B\text{-}\modd_{\ol{\cC}(\tti,\tti)}\text{-}B)}(N\circ_{A} A\circ_{A} M , (N\circ_{A} A\circ_{A} M)[k]) \simeq \HH_{\cC}^k(B),
\end{aligned}
\]
since $N\circ_{A} A\circ_{A} M \simeq B$.
\endproof

\begin{defn}
For a cyclic $2$-representation $\mbM$ of a fiat $2$-category $\cC$, we define $\HH^*_\cC(\mbM)$ to be $\HH^*_\cC(A)$ for $A$ an algebra $1$-morphism such that $\mbM\cong \mbM_A$.
\end{defn}

\subsection{The reduced Hochschild cohomology complex}
If $P_i=AF_iA$ for some $1$-morphism $F_i$ for all $i$ in the projective resolution of $A$, then by the free-forgetful adjunction for $A$-$A$-bimodules, we can construct $$\Hom_{\ol{\cC}(\tti,\tti)}(F_0,A)\to\Hom_{\ol{\cC}(\tti,\tti)}(F_1,A)\to\Hom_{\ol{\cC}(\tti,\tti)}(F_2,A)\to\dots$$ with the same cohomology.\par

We now give the isomorphisms explicitly.
Given $f\colon F_i\to A$, we have $$AF_iA\xrightarrow{\id_A\circh f\circh\id_A} AAA \xrightarrow{\mu\circv (\mu\circh\id_A)}A,$$ and given $g\colon AF_iA\to A$, we have $$F_i\xrightarrow{\iota\circh\id_{F_i} \circh\iota} AF_iA\xrightarrow{g}A.$$ Thus given a differential $d_i\colon AF_iA\to AF_{i-1}A$, and consequently the differential $$-\circ d_i\colon \Hom_{A\text{-}\modd_{\ol{\cC}(\tti,\tti)}\text{-}A}(AF_{i-1}A,A)\to \Hom_{A\text{-}\modd_{\ol{\cC}(\tti,\tti)}\text{-}A}(AF_{i}A,A),$$
 the corresponding differential  $\Hom_{\ol{\cC}(\tti,\tti)}(F_{i-1},A) \to \Hom_{\ol{\cC}(\tti,\tti)}(F_i,A)$ takes a morphism $f\colon F_{i-1}\to A$ to $$\mu\circv (\mu\circh \id_A)\circv (\id_A\circh f\circh\id_A)\circv d_i\circv (\iota\circh\id_{F_i}\circh\iota) \colon F_i\to A.$$

\subsection{Passing to cell-theoretic quotients}

Let $\cC$ be a quasi-fiat $2$-category and $\ceJ$ a fixed $2$-sided cell in $\cC$. Let $\cC/\cI_{\nleq\ceJ}$ be the quotient of $\cC$ by the $2$-ideal generated by all two-sided cells not less than or equal to $\ceJ$ and let $\cC_{\leq \ceJ}$ denote the $\ceJ$-simple quotient of $\cC$. Let further $\cC_{\ceJ}$ denote the $2$-full $2$-category on the same objects as $\cC_{\leq \ceJ}$, whose $1$-morphisms are those in the additive closure of the identities and $1$-morphisms in $\ceJ$. Then the natural $2$-functors in \cite[(2.11), (2.10), (2.8),(2.9)]{MMMTZbicat} provide biequivalences between the $2$-categories of finitary (resp.\ cyclic, transitive, simple) $2$-representations with apex $\ceJ$ of $\cC/\cI_{\nleq\ceJ}$  and of $\cC$.

By the dual for projective abelianisations of \cite[Theorem 4.26]{MMMTZbicat}, the $2$-category of cyclic $2$-representations of $\cC$ with exact morphisms is biequivalent to the bicategory $\cB\cB bimod_{\ol{\cC}}$ whose objects are algebra $1$-morphisms in $\ol{\cC}$ and whose morphism categories $\cB\cB bimod_{\ol{\cC}}(A,B)$ are the categories of biprojective $A$-$B$-bimodule $1$-morphisms.

Combining the previous two paragraphs with the observation that $A\text{-}\modd_{\ol{\cC}(\tti,\tti)}\text{-}A$ is the abelianisation of the category of projective objects in $\cB\cB bimod_{\ol{\cC}}(A,A)$, we obtain, for $A$ an algebra $1$-morphism in $\ol{\cC}(\tti,\tti)$ such that $\mbM_A$ has apex $\ceJ$, an equivalence of abelian categories
$$A\text{-}\modd_{\ol{\cC}(\tti,\tti)}\text{-}A \simeq A\text{-}\modd_{\ol{\cC/\cI_{\nleq \ceJ}}(\tti,\tti)}\text{-}A. $$

Moreover, by \cite[Theorem 4.22]{MMMTZbicat}, the category of simple $2$-representation of $\cC$ with apex $\ceJ$ is biequivalent to the category of simple $2$-representations of $\cC_{\leq \ceJ}$. If $\mbM_A$ is simple with apex $\ceJ$, then viewing this as a $2$-representation of $\cC_{\leq \ceJ}$, we have $A\in \cC_{\leq \ceJ}(\tti,\tti)$ by \cite[Theorem 4.19]{MMMTZbicat} and, again considering $A\text{-}\modd_{\ol{\cC}(\tti,\tti)}\text{-}A$ as the abelianisation of the category of projective objects in $\cB\cB bimod_{\ol{\cC}}(A,A)$ and using  \cite[Theorem 4.28, Remark 4.29]{MMMTZbicat}, this implies $$A\text{-}\modd_{\ol{\cC}(\tti,\tti)}\text{-}A\simeq A\text{-}\modd_{\ol{\cC_{\leq \ceJ}}(\tti,\tti)}\text{-}A\simeq A\text{-}\modd_{\ol{\cC_{\ceJ}}(\tti,\tti)}\text{-}A.$$

Putting these observations together, we have the following result.

\begin{psn}\label{passtoquot}
Let $\mbM_A$ be a finitary $2$-representation of $\cC$ with apex $\ceJ$. Then 
$$\HH^*_{\cC}(A) \cong \HH^*_{\cC/\cI_{\nleq\ceJ}}(A)$$
and if $\mbM_A$ is simple, then 
$$\HH^*_{\cC}(A) \cong \HH^*_{\cC_{\leq\ceJ}}(A)\cong \HH^*_{\cC_{\ceJ}}(A).$$
\end{psn}

\subsection{A replacement of the bar resolution}\label{sec:barres}

If $\mbM_A$ is a simple $2$-representation with apex $\ceJ$, we can pass to the $\ceJ$-simple quotient $\cC_{\leq \ceJ}$ by Proposition \ref{passtoquot} and without loss of generality assume that $\cC = \cC_{\leq \ceJ}$ and that, by \cite[Theorem 4.19]{MMMTZbicat}, the associated algebra $1$-morphism $A$ is in $\cC(\tti,\tti)$ , rather than the abelianisation. Then we can compute the Hochschild cohomology of $A$ using the usual bar resolution
\[
\cdots \to AAAA\to AAA\to AA\text{ of }A.
\] 
However, if $A \in \ol{\cC}$ is not in $\cC$ itself, then, while the term $AA$ is projective, the terms of the form $A^{\circ n}$, for $n \neq 2$ will not be projective, and so in this general case we need to construct a replacement for the bar resolution.

To this end, assume $A$ is given by $A_1\xrightarrow{a} A_0$ in the $\ol{\cC}(\tti,\tti)$. Then the multiplication map $\mu$ is given by $(\mu_0, (\mu_{01},\mu_{10}))$, i.e.
$$ 
\xymatrix{
A_0A_1\oplus A_1A_0 \ar_{ (\mu_{01},\mu_{10})}[d] \ar^{(\id_{A_0}\circh a, a\circh \id_{A_0})}[rrr]&&& A_0A_0\ar^{\mu_0}[d]\\
A_1  \ar^{a}[rrr] &&& A_0
}
$$

Note that, in particular $A$ is an $A_0$-$A_0$ bimodule with left, resp.\ right, actions given by
$\mu_{0\bullet} = (\mu_{01} , \mu_0)$, resp.\  $\mu_{\bullet 0} = (\mu_{10} , \mu_0)$, i.e.
$$ 
\xymatrix{
A_0A_1 \ar_{ \mu_{01}}[d] \ar^{\id_{A_0}\circh a}[rr]&& A_0A_0\ar^{\mu_0}[d]&&&
A_1A_0 \ar_{ \mu_{10}}[d] \ar^{ a\circh \id_{A_0}}[rr]&& A_0A_0\ar^{\mu_0}[d]\\  
A_1  \ar^{a}[rr] && A_0 &&&  A_1  \ar^{a}[rr] && A_0 
}
$$

\begin{lem}\label{barrep}
We have the beginning of a projective $A$-$A$-bimodule resolution of $A$ given by
$$\cdots \to AA_0A_0A\oplus AA_1A \xrightarrow{(\sigma, \id_A\circh a \circh \id_A)} AA_0A\xrightarrow{( \mu_{\bullet 0}\circh\id_A -  \id_A\circh\mu_{0\bullet})} AA $$
where $\sigma = \mu_{\bullet 0}\circh\id_{A_0A} - \id_A\circh\mu_0\circh\id_A + \id_{AA_0}\circh\mu_{0\bullet}$.
\end{lem}

\proof
Consider the usual bar resolution
$$\cdots \to AAAA\xrightarrow{\mu\circh\id_{AA}-\id_A\circh\mu\circh\id_A+\id_{AA}\circh\mu} AAA\xrightarrow{\mu\circh\id_A-\id_A\circh\mu}  AA$$ of $A$, which is exact, but fails to be projective in general.
However, using Lemma \ref{lem:ProjBimodDefn}, $AA = A\bbon_{\tti}A$ is a projective bimodule, a projective bimodule resolution of $AAA$ has begining $\cdots AA_1A\xrightarrow{\id_A\circh a\circh \id_A} AA_0A$ and the the first step in a projective bimodule resolution of $AAAA$ is given by $AA_0A_0A$. Splicing these projective bimodule resolutions of the components in the bar resolution together gives the desired result.
\endproof

\subsection{Hochschild cohomology in degree zero }

In this subsection, we give the analogue of the description of $\HH_{\cC}^0(A)$ as the \emph{centre} of an algebra.

\begin{psn}\label{centre}
For an algebra $1$-morphism $A = (A_1\xrightarrow{a}A_0)$ in $\ol{\cC}(\tti,\tti)$, there is an isomorphism of vector spaces
$$\HH_{\cC}^0(A) \cong \{f\in \Hom_{\ol{\cC}(\tti,\tti)}(\bbon_\tti, A)\,\mid \, \mu\circv (\id_A\circh f)=\mu\circv (f\circh \id_A)\}.$$
\end{psn}

\proof
Consider the projective bimodule resolution from Lemma \ref{barrep}. Then under the isomorphism 
$$
\Hom_{A\text{-}\modd_{\ol{\cC}(\tti,\tti)}\text{-}A}(AA,A)\to\Hom_{A\text{-}\modd_{\ol{\cC}(\tti,\tti)}\text{-}A}(AA_0A,A)\to \cdots .
$$ 
to 
$$
\Hom_{\ol{\cC}(\tti,\tti)}(\bbon_\tti,A)\xrightarrow{\delta_0}\Hom_{\ol{\cC}(\tti,\tti)}(A_0,A)\to\cdots 
$$
the $0$th cohomology consists of the kernel of $\delta_0$, where for $f\in \Hom_{\ol{\cC}(\tti,\tti)}(\bbon_\tti,A)$ $\delta_0(f)$ is given by the composition
\begin{equation}\label{redHHcom}
A_0 \xrightarrow{\iota\circh\id_{A_0}\circh \iota}AA_0A \xrightarrow{\mu_{\bullet 0}\circh\id_A -  \id_A\circh\mu_{0\bullet}} AA \xrightarrow{\id_A\circh f\circh \id_A}AAA \xrightarrow{\mu\circv(\mu\circh\id_A)}A.
\end{equation}

Denoting by $\pi$ the projection $\pi\colon A_0 \to A$, observe that  $$  \mu_{\bullet 0}\circv(\iota\circh\id_{A_0}) = \pi = \mu_{0\bullet}\circv(\id_{A_0}\circh \iota)$$
and hence
$$(\mu_{\bullet 0}\circh\id_A -  \id_A\circh\mu_{0\bullet})\circv(\iota\circh\id_{A_0}\circh \iota) = \pi\circh\iota - \iota\circh \pi$$
and $\delta_0(f)$ is given by the composition
$$A_0 \xrightarrow{\pi\circh\iota - \iota\circh \pi} AA \xrightarrow{\id_A\circh f\circh \id_A}AAA \xrightarrow{\mu\circv(\mu\circh\id_A)}A.$$

Using  the algebra axiom  $\mu\circv(\mu\circh\id_A) = \mu\circv(\id_A\circh\mu)$ and the interchange law,
we obtain
$$\mu\circv(\mu\circh\id_A)\circv(\id_A\circh f\circh \id_A)\circv (\pi\circh\iota) = \mu\circv(\id_A\circh\mu)\circv (\id_{AA}\circh\iota) \circv (\pi\circh f ) =\mu\circv(\pi\circh f )$$
and 
$$\mu\circv(\mu\circh\id_A)\circv(\id_A\circh f\circh \id_A)\circv(\iota \circh \pi) = \mu\circv(\mu\circh\id_A)\circv (\iota  \circh \id_{AA}) \circv ( f \circh \id_A) = \mu \circv ( f \circh \pi)$$
so $$\delta_0(f) = \mu\circv(\pi\circh f - f \circh \pi) =\mu\circv(\id_A\circh f - f \circh \id_A)\circv\pi. $$

Since $\pi$ is an epimorphism, $\delta_0(f)=0$  if and only if 
$$ \mu\circv(\id_A\circh f - f \circh \id_A) = 0,$$
from which the claim follows.
 \endproof
 
 A well-known result is that a simple algebra over an algebraically closed field has trivial centre. The next proposition provides the analogous statement for simple algebras (corresponding to simple $2$-representations) in fiat $2$-categories.
 
 \begin{psn} Let $A$ be an algebra $1$-morphism in $\cC(\tti,\tti)$, such that $\mbM_A$ is a simple $2$-representation. Then $\HH^0_\cC(A)\cong\bbk$.
\end{psn}
\proof
Let $\ceJ$ be the apex of $\mbM_A$. Then by Proposition \ref{passtoquot}, $\HH^0_\cC(A) \cong \HH^0_{\cC_{\leq\ceJ}}(A)$, so without loss of generality we may assume that $\cC$ is $\ceJ$-simple and, by \cite[Theorem 4.19]{MMMTZbicat}, that $A\in \cC(\tti,\tti)$ itself, rather than the abelianisation. 

Suppose $g\colon \bbon_\tti \to A$ gives rise to a nonzero element of $\HH^0_\cC(A)$, i.e.\ we have $\mu\circv (\id_A\circh g) = \mu\circv(g\circh \id_A)$.

Under the free forgetful adjunction, this corresponds to $f = \mu\circv(g\circh \id_A)$ in $\Hom_{\modd_{\ol{\cC}(\tti,\tti)}\text{-}A}(A,A)$. 

Since $f$ is annihilated by the differential, we have
\begin{align*}
0 & = \mu\circv(\id_A\circh f)\circv(\mu\circh \id_A-\id_A\circh\mu)\circv(\iota\circh\id_{AA})\\
& =\mu\circv(\id_A\circh f)\circv((\mu\circv ( \iota\circh\id_A))\circh \id_{A})
-\mu\circv(\iota\circh f)\circv\mu\\
& = \mu\circv(\id_A\circh f)-f\circv\mu
\end{align*}
and hence $\mu\circv(\id_A\circh f)=f\circv\mu$. Precomposing with $\id_A\circh \iota$, we derive that 
\begin{equation}\label{eq1}
f=\mu\circv(\id_A\circh f)\circv(\id_A\circh\iota).
\end{equation}

Assume $f\neq 0$, and consider the $\cC$-stable ideal in $\bproj_\cC(A)$ generated by $f$. Since $\bproj_\cC(A)$ is a simple $2$-representation, this ideal is all of $\bproj_\cC(A)$, and hence for any projective $A$-module $(X,\rho_X)$ in $\proj_{\cC(\tti,\ttj)}(A)$, there exists some $1$-morphism $F$ in $\cC$ such that $X$ is a direct summand of $FA$ with injection map $\sigma_X\colon X\to FA$ and projection map $\tau_X\colon FA\to X$ and such that $\tau_X\circv(\id_F\circh f)\circv\sigma_X=\lambda \id_X$ for some non-zero $\lambda\in\bbk$.
		
In more detail, simplicity implies that $\id_X=\sum_{i=1}^n a_i\circv(\id_{F_i}\circh f)\circv b_i$ for some $1$-morphisms $F_i$ and $2$-morphisms $a_i$, $b_i$ in $\cC$. We can rewrite this, using $F=\bigoplus_i^n F_i$, as $\id_X=\tau_X\circv(\id_F\circh f)\circv\sigma_X$ for some $2$-morphisms $\tau_X$ and $\sigma_X$.

Consider the diagram 
\[\begin{tikzcd}[ampersand replacement=\&]
	FA \& FAA \& FAA \& FA \\
	A \& AA \& AA \& A
	\arrow["{\on{id}_{FA} \circh \iota}", from=1-1, to=1-2]
	\arrow["{\on{id}_{F} \circh f}", shift left=3, curve={height=-18pt}, from=1-1, to=1-4]
	\arrow["{\on{id}_{FA} \circh f}", from=1-2, to=1-3]
	\arrow["{\on{id}_{F}\circh \mu}", from=1-3, to=1-4]
	\arrow["{\sigma_{A}}", from=2-1, to=1-1]
	\arrow["{\on{id}_{A} \circh\iota}"', from=2-1, to=2-2]
	\arrow["f", curve={height=18pt}, from=2-1, to=2-4]
	\arrow["{\sigma_{A} \circh \on{id}_{A}}", from=2-2, to=1-2]
	\arrow["{\on{id}_{A}\circh f}"', from=2-2, to=2-3]
	\arrow["{\sigma_{A} \circ \on{id}_{A}}", from=2-3, to=1-3]
	\arrow["\mu"', from=2-3, to=2-4]
	\arrow["{\sigma_{A}}", from=2-4, to=1-4]
\end{tikzcd}\]
Its top and bottom faces commute by \eqref{eq1}, the left-most two by bifunctoriality of $-\circh -$, and the rightmost middle face commutes since $\sigma_{A}$ is a morphism of modules. We find that
\[
\on{id}_{A} = \tau_{A} \circ (\on{id}_{F} \circh f) \circ \sigma_{A} = (\tau_{A} \circ \sigma_{A}) \circ \mu \circ (\on{id}_{A} \circh f) \circ (\on{id}_{A} \circh \iota) = \tau_{A} \circ \sigma_{A} \circ f
\]
and thus $f$ is invertible.

On the other hand, since $f= \mu\circv(g\circh \id_A) = \mu\circv (\id_A\circh g) $, $f$ is also a morphism of left $A$-modules, and hence a morphism of bimodules. Given indecomposability of $A$ as a bimodule over itself together with the fact that $\Bbbk$ is algebraically closed, we see that $f = \lambda \id_A +x$ for $\lambda\neq 0\in \Bbbk$ and some radical endomorphism $x$ of $A$. Since both $f$ and $\lambda \id_A$ are annihilated by the differential, so is $x$. On the other hand, any such $x$ must be invertible, which contradicts $x$ being in the radical. Thus $x=0$ and $f=\lambda \id_A$. The result follows.
\endproof

\subsection{Hochschild cohomology in degree one}

Here we provide an appropriate analogue of the description of first Hochschild cohomology as the space of derivations.

\begin{psn}\label{deriv}
For an algebra $1$-morphism $A = (A_1\xrightarrow{a}A_0)$ in $\ol{\cC}(\tti,\tti)$, there is an isomorphism of vector spaces
$$
\HH_{\cC}^1(A) \cong \{f\in \Hom_{\ol{\cC}(\tti,\tti)}(A, A)\,\mid \, f\circv\mu = \mu \circv(\id_A\circh f)+ \mu \circv(f\circh \id_A)\}/E
$$
where $$E = \{f\in \Hom_{\ol{\cC}(\tti,\tti)}(A, A) \,\vert\, \exists g\in \Hom_{\ol{\cC}(\tti,\tti)}(\bbon_\tti, A) \text{ with } f = \mu(\id_A\circh g-g\circh \id_A)\}.$$
\end{psn}

\proof
Consider the beginning of the projective $A$-$A$-bimodule resolution of $A$ given in Lemma \ref{barrep}.
Via the free-forgetful adjunction, $\HH_{\cC}^1(A)$ is given by the cohomology  in degree one of the reduced complex
$$\Hom_{\ol{\cC}(\tti,\tti)}(\bbon_\tti,A)\xrightarrow{\delta_0}\Hom_{\ol{\cC}(\tti,\tti)}(A_0,A)\xrightarrow{\delta_1} \Hom_{\ol{\cC}(\tti,\tti)}(A_0A_0\oplus A_1,A) $$
with $\delta_0$ given by the composition in \eqref{redHHcom} and, for $f\in \Hom_{\ol{\cC}(\tti,\tti)}(A_0, A)$, $\delta_1(f)$ defined to be the composition
$$ A_0A_0\oplus A_1\xrightarrow{\iota\circh\id_{A_0A_0\oplus A_1}\circh\iota}AA_0A_0A\oplus AA_1A\xrightarrow{(\sigma, \id_A\circh a \circh \id_A)}
AA_0A\xrightarrow{\id_A\circh f\circh \id_A}AAA\xrightarrow{\mu\circv(\mu\circh \id_A)}A.$$

For the restriction of $\delta_1(f)$ to $A_1$, we obtain the morphism $A_1\to A$ given by $f\circv a$ by the interchange law and the unitality axioms for $A$.

We next consider the summands of the restriction of $\delta_1(f)$ to $A_0A_0$ corresponding to the summands of $\sigma =(\mu_{\bullet 0}\circh\id_{A_0A} - \id_A\circh\mu_0\circh\id_A + \id_{AA_0}\circh\mu_{0\bullet})$ individually. 
Again denoting by  $\pi \colon A_0\to A$ the natural projection, we first observe that $$\mu_{\bullet 0}\circv(\iota\circh\id_{A_0}) = \pi = \mu_{0\bullet }\circv(\id_{A_0}\circh\iota).$$ 
This together with the interchange law, as well associativity and unitality of $A$,
 then yields
\begin{equation*}\begin{split}
\mu\circv(\mu\circh \id_A)&\circv(\id_A\circh f\circh \id_A)\circv(\mu_{\bullet 0}\circh\id_{A_0A})\circv(\iota\circh\id_{A_0A_0}\circh\iota)\\
&=\mu\circv(\mu\circh \id_A)\circv(\id_A\circh f\circh \id_A)\circv(\pi\circh\id_{A_0}\circh\iota)\\
&=\mu\circv( \id_A\circh\mu) \circv (\id_{AA}\circh \iota)   \circv(\pi\circh f)\\
&=\mu\circv(\pi\circh f),
\end{split}\end{equation*}

\begin{equation*}\begin{split}
\mu\circv(\mu\circh \id_A)&\circv(\id_A\circh f\circh \id_A)\circv(\id_A\circh\mu_0\circh\id_A )\circv(\iota\circh\id_{A_0A_0}\circh\iota)\\
&=\mu\circv(\mu\circh \id_A)\circv(\iota\circh\id_{A}\circh\iota)\circv f\circv\mu_0\\
&=f\circv\mu_0
\end{split}\end{equation*}

and

\begin{equation*}\begin{split}
(\mu\circv(\mu\circh \id_A))&\circv(\id_A\circh f\circh \id_A)\circv(\id_{AA_0}\circh\mu_{0\bullet})\circv(\iota\circh\id_{A_0A_0}\circh\iota)\\
&=\mu\circv(\mu\circh \id_A)\circv(\id_A\circh f\circh \id_A)\circv(\iota\circh\id_{A_0}\circh\pi)\\
&=\mu\circv( \mu\circh \id_A) \circv (\iota\circh \id_{AA})   \circv(f\circh \pi)\\
&=\mu\circv(f\circh \pi).
\end{split}\end{equation*}

Thus $\delta_1(f)=0$ if and only if both 
$$f\circv a = 0 \quad\text{ and } \quad f\circv\mu_0 = \mu\circv(\pi\circh f)+ \mu\circv(f\circh \pi).$$

The first condition implies that $f$ descends to an morphism $\bar{f}\colon A\to A$ given by
$$ 
\xymatrix{
A_1 \ar_{ 0}[d] \ar^{a}[rr]&& A_0\ar^{f}[d]\\  
A_1  \ar^{a}[rr] && A_0,
}
$$
such that $f=\bar{f}\circv \pi$. Inserting this into the second condition, and noting the equality $\pi\circv\mu_0= \mu\circ (\pi\circh \pi)$, we obtain
$$ \bar{f}\circv \pi\circv\mu_0 = \mu\circv(\id_A\circh\bar{f})\circv(\pi\circh \pi)+ \mu\circv(\bar{f}\circh\id_A)\circv(\pi\circh \pi).$$
Since $\pi$ is an epimorphism, it follows that $\bar{f}\colon A\to A$ gives rise to an element of $\HH^1_{\cC}(A)$ if and only if 
$ \bar{f}\circv\mu =  \mu\circv(\id_A\circh\bar{f})+\mu\circv(\bar{f}\circh\id_A)$, as claimed. The equivalence of two such $\bar{f},\bar{f'}$ provided their difference is in the subspace $E$ follows directly from the computation of $\delta_0$ in Proposition \ref{centre}.
\endproof

Even if the underlying field $\Bbbk$ is perfect, there are examples of simple algebra $1$-morphisms which are not separable, see Example \ref{VecGHH}.

\subsection{Hochschild cohomology in degree $2$.}

In order to compute $\HH^2$, we first extend the replacement of the bar resolution one step further. To this end, consider the kernel of the morphism $a$ in $\ol{\cC}(\tti,\tti)$, given by
$$\xymatrix{
A_3 \ar^{c}[rr]\ar[d] && A_2\ar^{b}[d] \\
0 \ar[rr]\ar[d] && A_1\ar^{a}[d] \\
0 \ar[rr]&& A_0.\\
}$$

Then the same proof as in Lemma \ref{barrep} (using the notation from there) shows that the next step in the replacement bar resolution is given by
$$\xymatrix{  \cdots \ar[r] &AA_0^{\circ 3}A\oplus AA_0A_1A  \oplus  AA_1A_0A\oplus AA_2A \ar^{\rho}[rr]  && AA_0A_0A\oplus AA_1A \ar^{(\sigma, \id_A \circh a \circh\id_A)}[d]\\ & AA&& \ar^{( \mu_{\bullet 0}\circh\id_A -  \id_A\circh\mu_{0\bullet})}[ll] AA_0A
}$$
where, omitting the symbols for horizontal composition in the matrix to save space,
$$\rho = \begin{pmatrix}\tau & \id_{AA_0} a  \id_A &  \id_{A}a  \id_{A_0A} & 0\\
0 & \id_{A}  \mu_{01} \id_A- \mu_{\bullet 0} \id_{A_1A}& \id_A\mu_{10} \id_{A} -\id_{AA_1}  \mu_{0\bullet}   & \id_A b  \id_A
\end{pmatrix}$$
and
$$\tau = \mu_{\bullet 0}\circh\id_{A_0A_0A} - \id_A\circh\mu_0\circh\id_{A_0A} + \id_{AA_0}\circh\mu_{0}\circh\id_A - \id_{AA_0A_0}\circh\mu_{0\bullet}.$$

\begin{psn}\label{HH2}
Let $\mbM=\mbM_A$ be a finitary $2$-representation of $\cC$ for some algebra $1$-morphism $A\in {\ol{\cC}(\tti,\tti)}$. Then $\HH^2_{\cC}(\mbM)$ is given by the quotient $C/B$ where 
$C$ is the subspace of $\Hom_{\ol{\cC}(\tti,\tti)}(A_0^{\circ 2}, A_0)\oplus \Hom_{\ol{\cC}(\tti,\tti)}(A_1, A_0)$ consisting of those $(g_0,g_1)$ such that 
\begin{align*}
0&=\pi\circv\left(\mu_0\circv(\id_{A_0}\circh g_0) -  g_0\circv(\mu_0\circh\id_{A_0})  +  g_0\circv(\id_{A_0}\circh \mu_0) - \mu_0\circv( g_0 \circh \id_{A_0})\right)\\
0&=\pi\circv\left( g_0\circv ( \id_{A_0}\circh a)+  g_1\circv \mu_{01} - \mu_0\circv( \id_{A_0}\circh  g_1)\right) \\
0&=\pi\circv\left( g_0\circv (a\circh \id_{A_0})+  g_1\circv \mu_{10} - \mu_0\circv( g_1\circh \id_{A_0})\right)\\
0&=\pi\circv  g_1\circv b
\end{align*}
and $B$ is the subspace consisting of elements $(g_0,g_1)$ such that 
\begin{align*}
 \pi\circv  g_0 &=\pi\circv \left( \mu_0\circv (\id_{A_0}\circh  f) -  f\circv \mu_0 + \mu_0\circv ( f\circh\id_{A_0})\right)\\
\pi\circv  g_1 &=  \pi\circv  f\circv a
\end{align*}
for some $f\in \Hom_{\ol{\cC}(\tti,\tti)}(A_0, A_0)$.
\end{psn}

\proof
Translating the replacement of the bar resolution above to the reduced Hochschild cohomology complex, we need to determine the cohomology in the middle term of 
$$ \Hom_{\ol{\cC}(\tti,\tti)}(A_0,A)\xrightarrow{\delta_1} \Hom_{\ol{\cC}(\tti,\tti)}(A_0A_0\oplus A_1,A) \xrightarrow{\delta_2}  \Hom_{\ol{\cC}(\tti,\tti)}(A_0^{\circ 3}\oplus A_0A_1  \oplus  A_1A_0\oplus A_2,A) $$
where, as before,
$$\delta_1(f) = \mu\circv(\mu\circh \id_A)\circv(\id_A\circh f\circh \id_A)\circv(\sigma, \id_A\circh a \circh \id_A)\circv(\iota\circh\id_{A_0A_0\oplus A_1}\circh\iota)$$
and $$\delta_2(g) =  \mu\circv(\mu\circh \id_A)\circv(\id_A\circh g\circh \id_A)\circv\rho\circv(\iota\circh\id_{A_0^{\circ 3}\oplus A_0A_1  \oplus  A_1A_0\oplus A_2}\circh\iota).$$

We first analyse under what conditions $\delta_2(g) = 0$. Write $g = (g_0,g_1)$ for $g_0\colon A_0^{\circ 2} \to A$, $g_1\colon A_1\to A$. Then $g\in \ker (\delta_2)$ if and only if the restriction of $\delta_2(g)$ to each summand in the domain is zero.

Consider first the restriction of $\delta_2(g)$ to $A_2$. This is given by
\begin{equation*}\begin{split} \mu\circv &(\mu\circh \id_A)\circv(\id_A\circh g_1\circh \id_A)\circv(\id_A\circh b\circh \id_A)\circv(\iota\circh\id_{A_2}\circh\iota)\\
&= \mu\circv(\mu\circh \id_A)\circv(\id_A\circh (g_1\circv b)\circh \id_A)\circv(\iota\circh\id_{ A_2}\circh\iota)\\
&= \mu\circv(\mu\circh \id_A)\circv(\iota\circh\id_{ A}\circh\iota)\circv(g_1\circv b)\\
&=g_1\circv b.
\end{split}\end{equation*}

Next, consider the restriction of $\delta_2(g)$ to $A_1A_0$. 
This is given by

\begin{equation*}\begin{split}
 \mu\circv &(\mu\circh \id_A)\circv(\id_A\circh (g_0, g_1)\circh \id_A)\circv\left( \begin{smallmatrix}\id_A\circh a\circh \id_{A_0A}\\  \id_A\circh\mu_{10} \circh\id_{A} -\id_{AA_1} \circh \mu_{0\bullet} \end{smallmatrix}\right)\circv(\iota\circh\id_{A_1A_0}\circh\iota)\\
&=\mu\circv(\mu\circh \id_A)\circv\left(\id_A\circh \left((g_0\circv (a\circh \id_{A_0})+ g_1\circv \mu_{10})\circh \id_A - g_1\circh \mu_{0\bullet}\right)\right)\circv(\iota\circh\id_{A_1A_0}\circh\iota)\\
&= g_0\circv (a\circh \id_{A_0})+ g_1\circv \mu_{10} - \mu\circv\left(g_1\circh \left( \mu_{0\bullet}\circv( \id_{A_0}\circh\iota)\right)\right)\\
&=  g_0\circv (a\circh \id_{A_0})+ g_1\circv \mu_{10} - \mu\circv(g_1\circh \pi)
\end{split}\end{equation*}
where the last step uses that  $\mu_{0\bullet}\circv( \id_{A_0}\circh\iota) = \pi$, the natural projection of $A_0$ to $A$ by construction of the maps in Section \ref{sec:barres}.

Similarly, we obtain that the restriction of $\delta_2(g)$ to $A_0A_1$ is given by
\begin{equation*}\begin{split}
\mu\circv &(\mu\circh \id_A)\circv(\id_A\circh (g_0, g_1)\circh \id_A)\circv\left(
 \begin{smallmatrix}\id_{AA_0}\circh a\circh \id_{A}\\   \id_{A} \circh \mu_{01} \circh\id_A- \mu_{\bullet 0}\circh \id_{A_1A}\end{smallmatrix}\right)\circv(\iota\circh\id_{A_0A_1}\circh\iota)\\
&=  g_0\circv ( \id_{A_0}\circh a)+ g_1\circv \mu_{01} - \mu\circv( \pi\circh g_1)
\end{split}\end{equation*}

Finally, consider the restriction of $\delta_2(g)$ to $A_0^{\circ 3}$. This is given by
\begin{equation*}\begin{split} \mu&\circv (\mu\circh \id_A)\circv(\id_A\circh g_0\circh \id_A)\\
&\circv(\mu_{\bullet 0}\circh\id_{A_0A_0A} - \id_A\circh\mu_0\circh\id_{A_0A} + \id_{AA_0}\circh\mu_{0}\circh\id_A - \id_{AA_0A_0}\circh\mu_{0\bullet})\\
&\circv(\iota\circh\id_{A_0^{\circ 3}}\circh\iota)\\
&=  \mu\circv (\mu\circh \id_A)\circv(\mu_{\bullet 0}\circh g_0\circh\id_{A})\circv(\iota\circh\id_{A_0^{\circ 3}}\circh\iota)\\
& \quad - g_0\circv(\mu_0\circh\id_{A_0}) + g_0\circv(\id_{A_0}\circh \mu_0) \\
& \quad - \mu\circv (\mu\circh \id_A)\circv(\id_{A}\circh g_0\circh \mu_{0\bullet})\circv(\iota\circh\id_{A_0^{\circ 3}}\circh\iota)\\
&= \mu\circv((\mu_{\bullet 0} \circv (\iota\circh \id_{A_0})\circh g_0) - g_0\circv(\mu_0\circh\id_{A_0}) \\
&\quad + g_0\circv(\id_{A_0}\circh \mu_0) - \mu\circv( g_0 \circh (\mu_{0\bullet} \circv (\id_{A_0}\circh \iota))) \\
&= \mu\circv(\pi\circh g_0) - g_0\circv(\mu_0\circh\id_{A_0})  + g_0\circv(\id_{A_0}\circh \mu_0) - \mu\circv( g_0 \circh \pi))
\end{split}\end{equation*}

Computing the image of $\delta_1$, the same computations as in Proposition \ref{deriv} show that a pair $(g_0,g_1)$ is in the image of $\delta_1$ if it is of the form
$$(g_0,g_1) = (\mu\circv (\pi\circh f) - f\circv \mu_0 + \mu\circv (f\circh \pi), f\circv a)$$
for some $ f\colon A_0\to A$.

Taking into account that any maps $g_0\colon A_0^{\circ 2} \to A$ and $g_1\colon A_1\to A$ necessarily factor over $A_0$ as
$g_0 = \pi\circv \hat g_0$ and $g_1 = \pi\circv \hat g_1$, the conditions above translate to 

\begin{equation*}\begin{split}
&\delta_2(g) = 0 \Leftrightarrow \\
&\begin{cases}
\pi\circv\left(\mu_0\circv(\id_{A_0}\circh \hat g_0) - \hat g_0\circv(\mu_0\circh\id_{A_0})  + \hat g_0\circv(\id_{A_0}\circh \mu_0) - \mu_0\circv(\hat g_0 \circh \id_{A_0})\right)=0 &\\
\pi\circv\left(\hat g_0\circv ( \id_{A_0}\circh a)+ \hat g_1\circv \mu_{01} - \mu_0\circv( \id_{A_0}\circh \hat g_1)\right) = 0&\\
\pi\circv\left(\hat g_0\circv (a\circh \id_{A_0})+ \hat g_1\circv \mu_{10} - \mu_0\circv(\hat g_1\circh \id_{A_0})\right)=0 &\\
\pi\circv \hat g_1\circv b =0&
\end{cases}
\end{split}\end{equation*}

Likewise 
\begin{equation*}\begin{split}
&(g_0,g_1)\in \im(\delta_1) \Leftrightarrow \\
&\begin{cases}
 \pi\circv \hat g_0 =\pi\circv \left( \mu_0\circv (\id_{A_0}\circh \hat f) - \hat f\circv \mu_0 + \mu_0\circv (\hat f\circh\id_{A_0})\right)&\\
\pi\circv \hat g_1 =  \pi\circv \hat f\circv a
\end{cases}
\end{split}\end{equation*}
for some $\hat f\colon A_0\to A_0$.

This completes the proof.
\endproof

\begin{rmk}\begin{enumerate}[(a)]\label{HH2rem}
\item\label{HH2rem1} If $\mbM=\mbM_A$ is simple, we can again without loss of generality assume that $\cC$ is $\ceJ$-simple with respect to the apex $\ceJ$ of $\mbM$ and that $A\in \cC(\tti,\tti)$, i.e. $A=A_0$ and $A_1, A_2=0$. Thus, the proposition simplifies to the usual description of second Hochschild cohomology as the subquotient of $\Hom_{{\cC}(\tti,\tti)}(AA, A) $ given by {\small
$$
\frac{ \{ g\,\vert\, \mu\circv(\id_{A}\circh g) -  g\circv(\mu\circh\id_{A})  +  g\circv(\id_{A}\circh \mu) - \mu\circv( g \circh \id_{A}=0\}  } {\{ g \,\vert\, \exists f  \in \Hom_{{\cC}(\tti,\tti)}(A, A) \text{ s.t. } g=\mu\circv (\id_{A}\circh  f) -  f\circv \mu + \mu\circv ( f\circh\id_{A}) \}.} 
$$}
\item\label{HH2rem2} The conditions only involving $g_1$ precisely imply that $g_1$ gives rise to an element in $\Ext^1_{\ol{\cC}(\tti,\tti)}(A,A)$. In particular, when all $g_0$ satisfying the first condition for being in a cocycle also satisfy the condition for being a coboundary, $\HH^2_{\cC}(A)$ is isomorphic to a subspace of $\Ext^1_{\ol{\cC}(\tti,\tti)}(A,A)$. An example where $\HH^2_{\cC}(A)\cong\Ext^1_{\ol{\cC}(\tti,\tti)}(A,A)$ is given in Section \ref{Mexample}.
\end{enumerate}
\end{rmk}

\subsection{$\HH^2$ and deformations}

By Remark \ref{HH2rem}\eqref{HH2rem1}, for a simple algebra $A$ or more generally, one that lives in $\cC(\tti,\tti)$, second Hochschild cohomology gives rise to infinitesimal deformations in precisely the way as in classical representation theory. If $A$ is not projective in $\ol{\cC}(\tti,\tti)$, we instead obtain a different construction, under the assumption that the lift $\mu_0$ of the multiplication map to $A_0$ induces an associative algebra structure on $A_0$.

Let $\Vect$ denote the monoidal category of finite-dimensional vector spaces. We define the bicategory $\cC \boxtimes \Vect$ by setting $\on{Ob}(\cC \boxtimes \Vect) = \on{Ob}(\cC)$, and, by setting $(\cC \boxtimes \Vect)(\tti, \ttj) := (\cC \boxtimes \Vect)(\tti, \ttj) \otimes_{\Bbbk} \Vect$. We will denote the $1$-morphism $(F,V)$ by $F \boxtimes V$ and similarly for $2$-morphisms, we define the horizontal composite $(G \boxtimes W)(F\boxtimes V) = GF \boxtimes (W \otimes_{\Bbbk} V)$, and similarly for $2$-morphisms. This is not a $2$-category, since $\Vect$ is not strict monoidal, and we lift the unitors and associators for $\cC \boxtimes \Vect$ in the obvious way from $\Vect$. Since the Hom-categories of $\cC$ are additive, the pseudofunctor from $\cC$ to $\cC \boxtimes \Vect$ sending $F$ to $F \boxtimes \Bbbk$, whose coherence $2$-morphisms are obtained from the unitors for $\Vect$, is a biequivalence. We fix a choice of a quasi-inverse to it, and, abusing notation, denote by $F \boxtimes V \in \cC$ the image of $F \boxtimes V \in \cC \boxtimes \Vect$ under this chosen biequivalence. In particular, $F \boxtimes V \simeq F^{\oplus \dim_{\Bbbk}V}$, and such an isomorphism is given by a choice of a basis for $V$; similarly, choosing bases we obtain an isomorphism $\Hom_{\cC(\tti,\ttj)}(F \boxtimes V, G\boxtimes W) \simeq \Hom_{\cC(\tti,\ttj)} (F,G)\otimes_{\Bbbk}\on{Mat}_{\dim_{\Bbbk}(W)\times \dim_{\Bbbk}(V)}(\Bbbk)$.

\begin{psn}\label{HH2def}
Assume that $A=(A_1\xrightarrow{a} A_0)$ is an algebra in $\ol{\cC}(\tti,\tti)$ with multiplication components $\mu_0,\mu_{01},\mu_{10}$ as in Section \ref{sec:barres}, such that $\mu_0$ defines an associative multiplication on $A_0$.

Assume further that $\boldsymbol{g} = (g_0,g_1)$ as in Proposition \ref{HH2} gives rise to an element of $\HH^2(A)$. Then there exist morphisms $$(h_{10},h_{01})\colon A_1A_0\oplus A_0A_1\to A_1$$ such that, letting $m$ denote multipliction in $\Bbbk[t]/(t^2)$, the diagram
\begin{equation}\label{defmult}
\xymatrix{
(A_1A_0\boxtimes  (\Bbbk[t]/(t^2))^{\otimes 2})\oplus (A_0A_1 \boxtimes  (\Bbbk[t]/(t^2))^{\otimes 2})\ar^{\alpha}[rr]\ar^{(\mu_{10}\boxtimes m ,\mu_{01}\boxtimes m )+(h_{10}\boxtimes tm,h_{01}\boxtimes tm)}[d]&& A_0A_0\boxtimes  (\Bbbk[t]/(t^2))^{\otimes 2}\ar^{\mu_0\boxtimes m +g_0\boxtimes tm}[d]\\
A_1\boxtimes \Bbbk[t]/(t^2) \ar^{a\boxtimes \id_{\Bbbk[t]/(t^2)} - g_1\boxtimes t}[rr] &&A_0\boxtimes \Bbbk[t]/(t^2)
}\end{equation}
where
$$\alpha = ((a\circh \id_{A_0})\boxtimes\id_{(\Bbbk[t]/(t^2))^{\otimes 2}},  (\id_{A_0}\circh a)\boxtimes \id_{(\Bbbk[t]/(t^2))^{\otimes 2}} )- ((g_1\circh \id_{A_0})\boxtimes (t\otimes \id_{\Bbbk[t]/(t^2)}) ,  (\id_{A_0}\circh g_1) \boxtimes(\id_{\Bbbk[t]/(t^2)}\otimes t) )$$
defines an associative multiplication, denoted by $\mu^{\boldsymbol{g}}$, on 
$$A^{\boldsymbol{g}}  = (A_1\boxtimes \Bbbk[t]/(t^2) \xrightarrow{a\boxtimes \id_{\Bbbk[t]/(t^2)} - g_1\boxtimes t} A_0\boxtimes \Bbbk[t]/(t^2)).$$
\end{psn}

\proof
The conditions
\begin{align*}
0&=\pi\circv\left( g_0\circv ( \id_{A_0}\circh a)+  g_1\circv \mu_{01} - \mu_0\circv( \id_{A_0}\circh  g_1)\right) \\
0&=\pi\circv\left( g_0\circv (a\circh \id_{A_0})+  g_1\circv \mu_{10} - \mu_0\circv( g_1\circh \id_{A_0})\right)\\
\end{align*}
imply that there exist $h_{01}, h_{10}$ such that
\begin{align*}
a\circv h_{01}&=\left( g_0\circv ( \id_{A_0}\circh a)+  g_1\circv \mu_{01} - \mu_0\circv( \id_{A_0}\circh  g_1)\right) \\
a\circv h_{10}&=\left( g_0\circv (a\circh \id_{A_0})+  g_1\circv \mu_{10} - \mu_0\circv( g_1\circh \id_{A_0})\right).\\
\end{align*}
With such $h_{01}, h_{10}$, we have
\begin{equation*}\begin{split}
(\mu_0&\boxtimes m +g_0\boxtimes tm)\circv ((a\circh \id_{A_0})\boxtimes \id_{(\Bbbk[t]/(t^2))^{\otimes 2}} - (g_1\circh \id_{A_0}) \boxtimes (t\otimes\id_{\Bbbk[t]/(t^2)})) \\
&=(\mu_0 \circv (a\circh \id_{A_0}))\boxtimes m +(g_0\circv (a\circh \id_{A_0}) -\mu_0\circv(g_1\circh \id_{A_0}))\boxtimes tm\\
&= (a\circv \mu_{10}) \boxtimes m + ( -g_1\circv\mu_{10}+ a\circv h_{10})\boxtimes tm\\
&= (a\boxtimes \id_{\Bbbk[t]/(t^2)}) \circv (\mu_{10}\boxtimes m) - (g_1\boxtimes t)\circv (\mu_{10}\boxtimes m) + (a\boxtimes \id_{\Bbbk[t]/(t^2)}) \circv(h_{10}\boxtimes tm)
\end{split}\end{equation*}
and similarly
\begin{equation*}\begin{split}
(\mu_0&\boxtimes m +g_0\boxtimes tm)\circv ((\id_{A_0}\circh a )\boxtimes \id_{(\Bbbk[t]/(t^2))^{\otimes 2}} - (\id_{A_0}\circh  g_1) \boxtimes (\id_{\Bbbk[t]/(t^2)} \otimes t)) \\
&=(\mu_0 \circv (\id_{A_0}\circh a ))\boxtimes m +(g_0\circv ((\id_{A_0}\circh a )) -\mu_0\circv(\id_{A_0}\circh  g_1))\boxtimes tm\\
&= (a\circv \mu_{01}) \boxtimes m + ( -g_1\circv\mu_{01}+ a\circv h_{01})\boxtimes tm\\
&= (a\boxtimes \id_{\Bbbk[t]/(t^2)}) \circv (\mu_{01}\boxtimes m) - (g_1\boxtimes t)\circv (\mu_{01}\boxtimes m) + (a\boxtimes \id_{\Bbbk[t]/(t^2)}) \circv(h_{01}\boxtimes tm)
\end{split}\end{equation*}
confirming the commutativity of \eqref{defmult}, so $\mu^{\boldsymbol{g}}$ is indeed a morphism from $A^{\boldsymbol{g}} A^{\boldsymbol{g}} $ to $A^{\boldsymbol{g}} $.

To check associativity, notice that computing $\mu^{\boldsymbol{g}}\circv(\mu^{\boldsymbol{g}}\circh \id_{A^{\boldsymbol{g}} }) - \mu^{\boldsymbol{g}}\circv(\id_{A^{\boldsymbol{g}} }\circh \mu^{\boldsymbol{g}})$
the component morphism $A_0^{\circ 3}\boxtimes \Bbbk[t]/(t^2)^{\otimes 3} \to A_0\boxtimes \Bbbk[t]/(t^2)$ is given by 
\begin{equation*}\begin{split}
(&\mu_0\circv(\mu_0 \circh \id_{A_0})) \boxtimes (m(m\otimes \id_{\Bbbk[t]/(t^2)} ))-
(\mu_0\circv( \id_{A_0}\circh \mu_0)) \boxtimes (m( \id_{\Bbbk[t]/(t^2)}\otimes m ))\\
&+(\mu_0\circv(g_0 \circh \id_{A_0})) \boxtimes (m(tm\otimes \id_{\Bbbk[t]/(t^2)} ))
+(g_0\circv(\mu_0\circh\id_{A_0})) \boxtimes (tm(m\otimes \id_{\Bbbk[t]/(t^2)} ))\\
&-(\mu_0\circv( \id_{A_0}\circh g_0)) \boxtimes (m( \id_{\Bbbk[t]/(t^2)}\otimes tm ))
-(g_0\circv( \id_{A_0}\circh \mu_0)) \boxtimes (tm( \id_{\Bbbk[t]/(t^2)}\otimes m ))\\
=&  (\mu_0\circv(\mu_0 \circh \id_{A_0} -  \id_{A_0}\circh \mu_0)) \boxtimes m(m\otimes \id_{\Bbbk[t]/(t^2)} )\\
&+ \big( \mu_0\circv(g_0 \circh \id_{A_0}) +g_0\circv(\mu_0\circh\id_{A_0}) - \mu_0\circv( \id_{A_0}\circh g_0)-g_0\circv( \id_{A_0}\circh \mu_0)  \big) \boxtimes tm(m\otimes \id_{\Bbbk[t]/(t^2)} )
\end{split}\end{equation*}
By our assumption that the $\mu_0$ defines an associative multiplication on $A_0$, 
$$\mu_0\circv(\mu_0 \circh \id_{A_0} -  \id_{A_0}\circh \mu_0) = 0. $$
The condition 
$$0=\pi\circv\left(\mu_0\circv(\id_{A_0}\circh g_0) -  g_0\circv(\mu_0\circh\id_{A_0})  +  g_0\circv(\id_{A_0}\circh \mu_0) - \mu_0\circv( g_0 \circh \id_{A_0})\right)$$
guarantees the existence of a morphism $\eta_1\colon A_0^{\circ 3}\to A_1$ such that 
$$\left(\mu_0\circv(\id_{A_0}\circh g_0) -  g_0\circv(\mu_0\circh\id_{A_0})  +  g_0\circv(\id_{A_0}\circh \mu_0) - \mu_0\circv( g_0 \circh \id_{A_0})\right) = a\circv \eta_1$$
which defines the required nullhomotopy for associativity.
\endproof

\begin{rmk}
Omitting the assumption that $\mu_0$ is associative, we obtain some conditions on the existence of a nullhomotopy for $\mu_0\circv(\mu_0 \circh \id_{A_0} -  \id_{A_0}\circh \mu_0) $, which composes to zero with some choice of representative of the element defined by $g_1 \in \Ext^1(A,A)$.
\end{rmk}

\section{Examples}\label{examplesec}

\subsection{Cell $2$-representations for strongly regular cells}\label{CellRigEx}

Recall the $2$-category $\cC_{B,X}$ of projective $B$-$B$-bimodules from \cite[Section~4.5]{MM3}, where $B$ is a finite-dimensional self-injective algebra and $X\subseteq Z(B)$ a subalgebra of the centre of $B$ which, under the isomorphism $Z(B) \cong \End_{B\text{-}\modd\text{-}B}(B)$, contains all elements which factor over $B\otimes_\Bbbk B$. More precisely, this has \begin{itemize}
\item an object $\tti$ for every connected component $B_i$ of $B$, which we identify with $B_i\text{-}\proj$;
\item $1$-morphisms in $\cC_{B,X}(\tti,\ttj)$ given by functors isomorphic to tensoring with a $B_i$-$B_j$-bimodule in the additive closure of $B_j\otimes_\Bbbk B_i$ if $\tti\neq \ttj$ and of $B_i\otimes_\Bbbk B_i \oplus B_i$ if $\tti= \ttj$,
\item all natural transformations of such functors as $2$-morphisms apart from the endomorphisms of $\bbon_\tti$, which correspond to the conponent of $X$ inside $Z(B_i)$.
\end{itemize}

\begin{thm}\label{rigidcell}
Let $\cC$ be a quasi-fiat $2$-category and $\ceL$ a left cell in a strongly regular two-sided cell. Then $\HH^n_\cC(\mbC_\ceL) = 0$ for $n>0$.
\end{thm}

\proof
By Proposition \ref{passtoquot}, we may without loss of generality assume that the two-sided cell $\ceJ$ containing $\ceL$ is the unique maximal cell of $\cC$, that $\cC$ is $\ceJ$-simple and that $\cC =\cC_{\ceJ}$, i.e. that its only $1$-morphisms are identities and those in $\ceJ$. Thus, by \cite[Theorem 32]{MM6}, we may assume that $\cC=\cC_{B,X}$ for some finite-dimensional self-injective algebra $B$ and a subalgebra $X\subseteq Z(B)$ of the centre of $B$. Moreover, an algebra $1$-morphism $A$ such that $\mbC_\ceL \cong \mbM_A$ can be chosen as $(eB)^*\otimes_\Bbbk eB$ for a primitive idempotent $e\in B$  with multiplication given by the natural evaluation map $$(eB)^*\otimes_\Bbbk eB\otimes_B(eB)^*\otimes_\Bbbk eB \to (eB)^*\otimes_\Bbbk eB$$ on the middle component, and the unit given by the unit $B\mapsto (eB)^*\otimes_\Bbbk eB$ of the adjunction. 
We claim that $A$ is projective as an $A$-$A$-bimodule. Since $AA = A\bbon A$ is a projective $A$-$A$-bimodule, it suffice to show that $A$ is a direct summand thereof as a bimodule. Let $e'$ be a primitive idempotent of $B$ such that $(eB)^*\cong Be'$. Then $AA \cong Be'\otimes_\Bbbk eBe'\otimes_\Bbbk eB \cong Be'\otimes_\Bbbk eB^{\oplus \dim eBe'}$ as an object of $\cC$. Moreover, since the left and right $A$-actions only depend on the outer tensor factors, this is indeed a decomposition of $A$-$A$-bimodules, proving the claim. Thus $A$ is separable and $\HH^n_\cC(\mbC_\ceL) = \HH^n_\cC(A) = 0$ for $n>0$.
\endproof

\subsection{Inflations of Cell $2$-Representations}\label{CellInfEx}

We consider the algebra viewpoint of the inflations considered in \cite{MM6}. Let $\mbM = \mbM_A$ be a cyclic $2$-representation of $\cC$ with associated algebra $1$-morphism $(A,\mu_A,\iota_A)$, and let $(R,\mu_R, \iota_R)$ be a finite-dimensional $\Bbbk$-algebra with multiplication $\mu_R$ and unit (homomorphism) $\iota_R$. 

Recall the inflation of $\mbM$ by $R$, denoted by $\mbM^{\boxtimes R}$, which is defined by setting 
\begin{itemize}
\item $\mbM^{\boxtimes R}(\tti)=\mbM(\tti)\boxtimes R\text{-}\proj$ for an object $\tti\in \cC$;
\item $\mbM^{\boxtimes R}(F)=\mbM(F)\boxtimes \Id_{R\text{-}\proj}$ for a $1$-morphisms $F$ in $\cC$;
\item $\mbM^{\boxtimes R}(\alpha)=\mbM(\alpha)\boxtimes  \id_{\Id_{R\text{-}\proj}}$ for a $2$-morphisms $\alpha$ in $\cC$.
\end{itemize}

Choose a basis $\{r_i \,\vert\, i=1,\ldots,  \dim_\Bbbk R\}$ of $R$ such that $r_1 = 1_R$ with structure constants $r_ir_j = \sum_{k\in I} c_{ij}^kr_k$. Denote by $\rho_{r_i}\in\Hom_{R\text{-}\proj}(R, R)$ the map given by right multiplication with $r_i$. 

By definition of $A\boxtimes R$ in $\ol{\cC}\boxtimes \Vect$, this object clearly admits an algebra $1$-morphism structure in $\ol{\cC}\boxtimes \Vect$, and hence also in $\ol{\cC}$. Indeed, denoting the multiplication in $A$ by $\mu_{A}$ and the multiplication in $R$ by $\mu_{R}$, the multiplication in this structure is given by $\mu_{A} \boxtimes \mu_{R}$. We denote this algebra $1$-morphism in $\ol{\cC}$ by $A^{\boxtimes R}$.

Fix a decomposition $A^{\boxtimes R} = \bigoplus_{i=1}^{\dim_\Bbbk R} A_{(i)}$ and denote by $p_i\colon A^{\boxtimes R} \leftrightarrows A \, :j_i$ the corresponding projection respectively injection. Then multiplication on $A^{\boxtimes R}$ is given by
$\sum   c_{ij}^k j_k   \circv\mu_A\circv(p_i\circh p_j)$ and the unit is given by $j_1\circ \iota_A$.

\begin{lem}
The $2$-representation $\mbM^{\boxtimes R}$ is equivalent to $\mbM_{A^{\boxtimes R}}$.
\end{lem}

\proof
The $2$-representation $\mbM^{\boxtimes R}$  is generated by the object $A\boxtimes R \in \mbM(\tti)\boxtimes R\text{-}\proj$, so we wish to compute $B:= {}_{\mbM^{\boxtimes R}}[A\boxtimes R,A\boxtimes R  ]$. One directly checks that for $F\in \cC(\tti,\tti)$, 
\begin{equation*}\begin{split}
\Hom_{\ol{\cC}(\tti,\tti)}(F, {}_{\mbM^{\boxtimes R}}[A\boxtimes R,A\boxtimes R  ])&\cong
 \Hom_{\mbM^{\boxtimes R}(\tti)}(\mbM^{\boxtimes R}(F)(A\boxtimes R), A\boxtimes R)\\
 &=\Hom_{\mbM^{\boxtimes R}(\tti)}(FA\boxtimes R, A\boxtimes R) \\
 &= \Hom_{\mbM(\tti)}(FA, A)\otimes_\Bbbk \Hom_{R\text{-}\proj}(R, R)\\
&  \cong \Hom_{\ol{\cC}(\tti,\tti)}(F,  {}_{\mbM}[A,A])\otimes_\Bbbk  \Hom_{R\text{-}\proj}(R, R) \\
  &\cong \Hom_{\ol{\cC}(\tti,\tti)}(F, A)\otimes_\Bbbk  \Hom_{R\text{-}\proj}(R, R)  \\
  &\cong \Hom_{\ol{\cC}(\tti,\tti)}(F, A)\otimes_\Bbbk  (\bigoplus_{i=1}^{\dim_\Bbbk R} \Bbbk\rho_{r_i})\\
 &\cong \Hom_{\ol{\cC}(\tti,\tti)}(F, A^{ \oplus \dim_\Bbbk R})
\end{split}\end{equation*}
hence, as an object of $\ol{\cC}(\tti,\tti)$, we indeed have $B\cong A^{\boxtimes R}$.

To check that the multiplication morphism is given as defined above, consider the isomorphism
$$\Hom_{\ol{\cC}(\tti,\tti)}(A, B) \cong \Hom_{\ol{\cC}(\tti,\tti)}(A, A)\otimes_\Bbbk \Hom_{R\text{-}\proj}(R, R) $$
and define $j_i \colon A \to B$ as the morphism corresponding to $\id_A \otimes \rho_{r_i}$, i.e. the embedding into the $i$th component of the decomposition determined by the choice of basis of $R$. Note that, in particular, the $j_i$ are split mono, and for each $i$, let $p_i$ denote a splitting, such that $p_ij_i=\id_A$, $p_ij_l=0$ for $i\neq l$, and $\id_B = \sum_i j_ip_i$.

Then, for any $F\in {\cC}(\tti,\tti)$, the isomorphism above is explicitly given by
\begin{align*}
\Hom_{\ol{\cC}(\tti,\tti)}(F, B) &\cong \Hom_{\ol{\cC}(\tti,\tti)}(F, A)\otimes_\Bbbk \Hom_{R\text{-}\proj}(R, R) \\
j_i\circv \phi &\leftrightarrow \phi\otimes  \rho_{r_i}.
\end{align*}
This immediately yields that the unit morphism is given by $j_1\circv \iota_A$.

Considering the commutative diagram
$$\xymatrix{
\Hom_{\ol{\cC}(\tti,\tti)}(A, B) \ar@/^/^{-\circv p_i}[rr]\ar@{<->}^{\sim}[d]&&\Hom_{\ol{\cC}(\tti,\tti)}(B, B) \ar@/^/^{-\circv j_i}[ll]\ar@{<->}^{\sim}[d]\\
\Hom_{\ol{\cC}(\tti,\tti)}(A, A)\otimes_\Bbbk \Hom_{R\text{-}\proj}(R, R) \ar@{<->}^{\sim}[d]&&  \Hom_{\mbM^{\boxtimes R}(\tti)}(\mbM^{\boxtimes R}(B)(A\boxtimes R), A\boxtimes R) \ar@{<->}^{\sim}[d]\\
\Hom_{\mbM(\tti)}(\mbM(A)A, A)\otimes_\Bbbk \Hom_{R\text{-}\proj}(R, R) \ar@/^/^{-\circv (\mbM(p_i)_A\otimes \id_R)}[rr]&& \Hom_{\mbM(\tti)}(\mbM(B)A, A)\otimes_\Bbbk \Hom_{R\text{-}\proj}(R, R)\ar@/^/^{-\circv (\mbM(j_i)_A\otimes \id_R)}[ll]
}$$
we see that the evaluation $\ev_B\colon\mbM^{\boxtimes R}(B)(A\boxtimes R)\to A\boxtimes R$ corresponds to
$$\sum_{i=1}^{\dim_\Bbbk R}(\ev_A\circv \mbM(p_i )_A)\otimes \rho_{r_i}$$
under the isomorphism
$$ \Hom_{\mbM^{\boxtimes R}(\tti)}(\mbM^{\boxtimes R}(B)(A\boxtimes R), A\boxtimes R) \cong  \Hom_{\mbM(\tti)}(\mbM(B)A, A)\otimes_\Bbbk \Hom_{R\text{-}\proj}(R, R)$$
Thus, under the isomorphism 
$$ \Hom_{\mbM^{\boxtimes R}(\tti)}(\mbM^{\boxtimes R}(BB)(A\boxtimes R), A\boxtimes R) \cong  \Hom_{\mbM(\tti)}(\mbM(BB)A, A)\otimes_\Bbbk \Hom_{R\text{-}\proj}(R, R)$$
the morphism $\ev_B\circv(\id_B\circh\ev_B)$ corresponds to
$$\sum_{i,l=1}^{\dim_\Bbbk R}(\ev_A\circv(\id_A\circh \ev_A)\circv \mbM(p_i \circh p_l)_A)\otimes \rho_{r_ir_l} =\sum_{i,l, k=1}^{\dim_\Bbbk R} c_{il}^k(\ev_A\circv(\id_A\circh \ev_A)\circv \mbM(p_i \circh p_l)_A)\otimes \rho_{r_k}  .$$
This, in turn, corresponds to
$$\sum_{i,l, k=1}^{\dim_\Bbbk R} c_{il}^k(\mu_A\circv (p_i \circh p_l))\otimes \rho_{r_k} $$
under the isomorphism
$$\Hom_{\mbM(\tti)}(\mbM(BB)A, A)\otimes_\Bbbk \Hom_{R\text{-}\proj}(R, R) \cong \Hom_{\ol{\cC}(\tti,\tti)}(BB, A)\otimes_\Bbbk \Hom_{R\text{-}\proj}(R, R) $$
which, finally under the isomorphism $$\Hom_{\ol{\cC}(\tti,\tti)}(BB, A)\otimes_\Bbbk \Hom_{R\text{-}\proj}(R, R) \cong \Hom_{\ol{\cC}(\tti,\tti)}(BB, B)$$
corresponds to 
$$\mu_B =  \sum_{i,l, k=1}^{\dim_\Bbbk R} c_{il}^k j_k(\mu_A\circv (p_i \circh p_l))$$
as claimed.
\endproof

For simple $2$-representations $\mbM$, we can compute the Hochschild cohomology of $\mbM^{\boxtimes R}$ as a tensor product.

\begin{psn}
Let $\mbM = \mbM_A$ be a simple $2$-representation of $\cC$ and $R$ a finite-dimensional $\Bbbk$-algebra. Then
$$\HH^*_{\cC}(\mbM^{\boxtimes R})\cong  \HH^*_{\cC}(\mbM)\otimes_\Bbbk \HH^*( R).$$
\end{psn}

\proof
Due to simplicity of $\mbM$, we can compute $\HH^*_{\cC}(\mbM) = \HH^*_{\cC}(A)$ via the usual bar resolution of $A$. Then the same proof as in \cite[Lemma 3.1]{LZ} shows that the usual bar resolution of $A\boxtimes R$, whose $i$th component is $(A\boxtimes R)^{\circ i+2}\cong A^{\circ i+2} \boxtimes R^{\otimes i+2}$, is homotopy equivalent to the tensor product of the bar resolutions, whose $i$th component is given by $\bigoplus_{j=0}^i A^{\circ j+2}\boxtimes R^{\otimes (i+2-j)}$, implying the result.
\endproof

Using Theorem \ref{rigidcell}, we immediately obtain the following corollary, which, using \cite[Theorem 4]{MM6}, completely describes the Hochschild cohomology of isotypic transitive $2$-representations for $\ceJ$-simple $\cC$ with strongly regular apex $\ceJ$.

\begin{cor}
If $\mbM$ is a cell $2$-representation for a left cell inside a strongly regular two-sided cell, then $\HH^*_{\cC}(\mbM^{\boxtimes R})\cong \HH^*( R)$.
\end{cor}

\subsection{$\mbC_{\{\bbon\}}$ extended by $\mbC_{\{F\}}$ for $\cC_D$}\label{Mexample}
Let $D=\bbk[x]/(x^2)$and consider $\cC_D$. Set $F$ be the indecomposable $1$-morphism corresponding to tensoring with  $D\otimes_\Bbbk D$. The indecomposable $1$-morphisms in $\cC_D$ are $\bbon$ and $F$ and each forms a left, right and  two-sided cell. Let $\mbM$ be the extension between the two cell $2$-representations considered in \cite[Subsection 4.1]{CM}. In loc.\ cit.\ we computed the corresponding coalgebra $1$-morphism; the dual computation shows that the algebra $1$-morphism $A$ in the projective abelianisation is given by the extension of the simple top $L_\bbon$ of $\bbon$ by the simple top $L_F$ of $F$, or equivalently, as the object $\bbon\xrightarrow{a}\bbon$ in the abeliansation where $a$ corresponds to multiplication by $x$ when viewed as an endomorphism of the identity bimodule $D$. Multiplication is given by the diagram
$$ 
\xymatrix{
\bbon\bbon\oplus \bbon\bbon\ar_{ (\id_{\bbon},\id_{\bbon})}[d] \ar^{(\id_{\bbon}\circh a, a\circh \id_{\bbon})}[rrr]&&& \bbon\bbon\ar^{\id_{\bbon}}[d]\\
\bbon  \ar^{a}[rrr] &&& \bbon
}
$$
One easily checks that in the extended bar resolution, the morphism $b$ is again the same as $a$.

Given that $\Hom_{\ol{\cC_D}(\bullet,\bullet)}(A,A) = \Bbbk\cdot \id_{A}$, it is easy to see that there are no derivations.

Since the only morphisms $\bbon\bbon = \bbon \to A$ are given by $\lambda\pi= \pi\circ (\lambda \id_\bbon)$ for $\lambda\in \Bbbk$, it is easy to check that the conditions given for $\delta_2(g_0,g_1)=0$ in Proposition \ref{HH2} hold for any pair $(g_0,g_1)=(\lambda_0 \pi, \lambda_1 \pi)$. On the other hand, such $(g_0,g_1)$ is in the image of $\delta_1$ only provided that there exists an $f = \lambda \pi$ such that $(\lambda_0 \pi, \lambda_1 \pi) = (\lambda \pi, 0)$. Thus, the second Hochschild cohomology is given by all morphisms from $A_1$ to $A$. These trivially annihilate $b$ and do not factor over $a$, and hence
$$\HH^2_{\cC_D}(\mbM) \cong \Ext^1_{\ol{\cC_D}(\bullet,\bullet)}(A,A).$$
Moreover, considering the projective resolution of $A$ in $\ol{\cC_D}(\bullet,\bullet)$ given by 
$$\cdots \bbon\xrightarrow{a}\bbon\xrightarrow{a}\bbon\xrightarrow{a}\bbon,$$ 
it is easy to see that $\Ext^1_{\ol{\cC_D}(\bullet,\bullet)}(A,A)\cong \Bbbk$.

\subsection{Nontrivial Hochschild cohomology for simple algebra $1$-morphisms in ${\cV}ec_G$ }\label{VecGHH}

In this subsection, we assume that the characteristic of $\Bbbk$ is $p>0$ and let $G=C_p = \langle g \,\vert\, g^p=1\rangle$ be the cyclic group of order $p$.  Consider the category ${\cV}ec_G$ of finite dimensional $G$-graded $\Bbbk$-vector spaces, denoting the indecomposable $1$-morphisms by $F_{g^a}$.

Consider the algebra $1$-morphism $A= \bigoplus_{g\in G}F_g$, which gives rise to the trivial (rank one) $2$-representation and is clearly simple. 

We claim that this has nontrivial first Hochschild cohomology. Indeed, first notice that any morphism $f$ in $ \Hom_{\ol{\cC}(\tti,\tti)}(A, A)$ is simply a collection of scalars $(f_0,\ldots,f_{p-1})$ where the component of $f$ from $F_{g^i}$ to $F_{g^j}$ is zero if $i\neq j$ and given by $f_i \id_{F_{g^i}}$ if $i=j$. For $f$ to be a derivation, we require $f_0=0$ and $f_a = a f_1$ for all $a = 2,\ldots, p$. On the other hand, it is easy to see that for any $f_1\in \Bbbk$, the morphism defined by $(0,f_1,2f_1, \ldots, (p-1)f_1)$ is indeed a derivation and that, moreover, there are no inner derivations. Therefore, any such $f_1$ gives rise to a distinct element of $\HH^1(A)$ and $\HH^1(A) \cong \Bbbk.$

It is easy to see that similar reasoning implies nontrivial first Hochschild cohomology for bigger groups $G$ and algebra $1$-morphisms given by more general subgroups of order divisible by $p$.

\end{document}